\documentclass[11pt]{article}

\usepackage{amsmath,amsthm,epsfig,amssymb}

\newtheorem{proposition}{Proposition}
\newtheorem{theorem}{Theorem}
\newtheorem{lemma}{Lemma}
\newtheorem{corollary}{Corollary}
\newtheorem{remark}{Remark}
\newtheorem{example}{Example}

\newtheorem{definition}{Definition}

\newenvironment{pf}{{\bf Proof:\
}}{\raisebox{-4pt}{\rule{6pt}{6pt}}\vspace{4pt plus 8pt minus 1pt}\vskip
8pt}

\newenvironment{prooft}{{\leftmargin=20pt\rightmargin=0pt} {\bf Proof of
Theorem:\ }}{\vspace{4pt plus 2pt minus 1pt}}

\newcommand{\bth}{\begin{theorem}}
\newcommand{\bpr}{\begin{proposition}}
\newcommand{\epr}{\end{proposition}}
\newcommand{\bco}{\begin{corollary}}
\newcommand{\eco}{\end{corollary}}
\newcommand{\ble}{\begin{lemma}}
\newcommand{\ele}{\end{lemma}}
\newcommand{\bre}{\begin{remark}\rm}
\newcommand{\ere}{\end{remark}}
\newcommand{\bex}{\begin{example}\rm}
\newcommand{\eex}{\end{example}}
\newcommand{\bres}{{\bf Remarks}\begin{enumerate}}
\newcommand{\eres}{\end{enumerate}}
\newcommand{\bexs}{{\bf Examples}\begin{enumerate}\begin{enumerate}}
\newcommand{\eexs}{\end{enumerate}\end{enumerate}}
\newcommand{\bde}{\begin{definition}\rm}
\newcommand{\ede}{\end{definition}}
\newcommand{\bpf}{\begin{pf}}
\newcommand{\epf}{\end{pf}}
\newcommand{\bproof}{\begin{pf}}
\newcommand{\eproof}{\end{pf}}
\newcommand{\bpft}{\begin{prooft}\begin{itemize}}
\newcommand{\epft}{\end{itemize}\end{prooft}}
\def\aut{\hbox{Aut}}

\def\la#1{\hbox to #1pc{\leftarrowfill}}
\def\ra#1{\hbox to #1pc{\rightarrowfill}}
\def\fract#1#2{\raise4pt\hbox{$ #1 \atop #2 $}}

\def\decdnar#1{\phantom{\hbox{$\scriptstyle{#1}$}}
\left\downarrow\vbox{\vskip15pt\hbox{$\scriptstyle{#1}$}}\right.}

\def\za{\vrule height6pt width4pt depth1pt}
\def\lrar{{\ra 2}}

\def\bbc{\mathbb C}
\def\bbp{\mathbb P}

\def\bbz{\mathbb Z}
\def\bbr{\mathbb R} 
\def\bbq{\mathbb Q}
 
\def\bbf{\mathbb F} 

\input xy
\xyoption{all}
\begin{document}
\bibliographystyle{plain}


\font\twelverm=cmr10 at 12pt
\font\svntnrm=cmr21
\centerline{\svntnrm On the Group of Automorphisms of}
\medskip
\centerline{\svntnrm Cyclic Covers of the Riemann Sphere}
\renewcommand{\baselinestretch}{1.0}

\vskip 20pt
\centerline{\large Sadok Kallel~~and~~Denis Sjerve
\footnote{Research supported by NSERC grant A7218}}

\begin{abstract}
This note uses some recent calculations of Conder and Bujalance (on
classifying finite index group extensions of Fuchsian groups with
abelian quotient and torsion free kernel) in order to determine the
full automorphism groups of some cylic coverings of the line
(including curves of Fermat and Lefschetz type). The answer is complete
for cyclic covers that branch over three points.
\end{abstract}


{\bf\Large\S1 Introduction}

Let $C$ be a Riemann surface of genus $g\geq 2$. As is known, the
group of automorphisms of such a curve is finite. If the curve is
given by an explicit affine equation in $\bbc^2$, a very compelling
problem is to determine the symmetry group from the equation.  In
this paper we answer this problem completely for certain cyclic covers
of the line.

A cyclic covering of the line is a curve $C$ with an affine equation in
$\bbc^2$ given by
$$C: y^n = (x-e_1)^{a_1}(x-e_2)^{a_2}\cdots (x-e_k)^{a_k}$$ 
We only deal with {\sl irreducible} curves (putting a small
restriction on $n$ and the $a_i$'s). The special cases when $n=2$ or
$n=p$ (an odd prime) are known as hyperelliptic curves or
$p$-elliptic curves (see \cite{KSS}). In \cite{Na}, Namba conjectured that the
isomorphism type of a cyclic cover $y^n=f(x)$ is determined entirely by
the branching of the projection $(x,y)\in C\mapsto x$.  This has
been answered affirmatively for $n=p$ by Nakajo \cite{Nk}.  A
general discussion of these curves and their properties is given in
\S3.  A consequence is that Namba's conjecture is true as far as the
automorphism group of $C$ is concerned, that is the isomorphism type of
the group of automorphisms of $C$ is determined by the branching.

A cyclic cover $C$ is clearly endowed with an action by the cyclic
group $\bbz_n$, $(x,y)\mapsto (x,\zeta y)$, where $\zeta$ is a
primitive $n$-th root of unity. If $\aut (C)$ is the group of
symmetries of $C$ then $\bbz_n\subset \aut (C)$. Our aim then becomes
to determine all possible extensions of $\bbz_n$ that can occur as
automorphism groups of cyclic galois covers.  To this end, we rely on
techniques in Fuchsian group theory and on pivotal recent calculations
of Bujalance, Conder and Cirre (\cite{BC}, \cite{BCC}).

The main result of this paper classifies all group actions on cyclic
covers which ramify over three distinct points of $\bbp^1$ (the choice
of these points is immaterial since PSL$_2(\bbc )$ is
$3$-transitive). These covers are given by
$$C: y^n = x^a(x-1)^b(x+1)^c,\ \hbox{where}\ 1\leq a,b,c\leq n-1, \
a+b+c\equiv 0\ (mod\ n)$$ 
We refer to them as cyclic {\sl belyi} covers. A curve given as above
is irreducible if $[n,a,b,c]=1$ and its genus is 
$\displaystyle g=(2+n-[n,a]-[n,b]-[n,c])/2$, see \S3.  We use the 
notation $[r,s]$ for the greatest common divisor of $r,s.$

Given $n>3$, we say that $(a,b,c)$ is equivalent to $(a',b',c')$ if
there is $k$ prime to $n$ and a permutation $\tau$ on $3$ letters such
that $a'\equiv k\tau(a)\ (mod\ n)$, $b'\equiv k\tau(b)\ (mod\ n)$
and\linebreak $c'\equiv k\tau(c)\ (mod\ n)$.  It is well-known
(Nielsen) that any two equivalent triples yield isomorphic cyclic
covers. So as far as automorphism groups are concerned, only
equivalence classes of triples will matter. 
We can now state our main result.

\bth \label{maintheorem} 
Let $C: y^n = x^a(x-1)^b(x+1)^c$ be a cyclic $\bbz_n$ Galois cover of
the line with $a+b+c\equiv 0\ (mod\ n)$, $1\le  a,b,c\le n-1$ and
$[n,a,b,c]=1$.  We fix $n\geq 4$ and let $G$ be the full
automorphism group of $C$. Then $G$ is completely determined by the
values of $n$, $a,$ $ b$ and $c$ as follows.\\
(i) For the following equivalence classes of triples $(a,b,c)$, 
$G$ is given by
{\small
$$\begin{tabular}{|c|l|c|c|c|c|c|c|}\hline
&n&a&b&c&g&$|G|$&G\\
\hline
\hbox{A.1}&$\hbox{odd}$&$1$&$1$&$n-2$&$\frac{n-1}{ 2}$&$2n$&
$\bbz_{2n}$\\
\hline
\hbox{A.2}&$\hbox{even}$&$1$&$1$&$n-2$&$\frac{n}{ 2}-1$&$4n$&
(central $\bbz_2): D_{2n}$\\
\hline
B.1&${{n\not\equiv 0 (8)}\atop {n\ne 12}}$&$1$&${b\neq 1, b^2\equiv 1 (n)}$
&$n-1-b$&$\frac{1}{ 2}(n-[n,b+1])$&$2n$&$\bbz_n\ltimes\bbz_2$\\
\hline
B.2&$8|n, n > 8$&$1$&${\frac{n}{2}-2}$
&$\frac{n}{2}+1$&$\frac{n}{2}-1$&$4n$&$(\bbz_n\ltimes\bbz_{2}):\bbz_2$\\
\hline
B.3&$8$&$1$&$2$&$5$&$3$&$96$&$(\bbz_4\oplus\bbz_4)\ltimes S_3$\\
\hline
C.1&${{n > 7,~n\equiv 1(2)}\atop {\exists\ p|n, p\equiv 1(3)}}$
&$1$&$b\neq 1, [b,n]=1$&$b^2$&
$\frac{n-1}{ 2}$&$3n$&$\bbz_n\ltimes\bbz_3$\\
\hline
C.2&$7$&$1$&$2$&$4$&$3$&$168$&PSL$(2,7)$\\
\hline
D.1&$12$&$1$&$3$&$8$&$3$&$48$& (central $\bbz_4 ): A_4$\\
\hline
E.1&$8$&$1$&$3$&$4$&$2$&$48$&GL$(2,3)$\\
\hline
E.2&$12$&$1$&$4$&$7$&$4$&$72$&(central $\bbz_3):S_4$ \\
\hline
E.3&$24$&$1$&$4$&$19$&$10$&$144$& (central $\bbz_6):S_4$\\
\hline
\end{tabular}$$
}

We use the notation $H\ltimes K$ to denote a semi-direct product of $H$ by $K$ 
and $H:K$ to denote a non-split extension of $H$ by $K$. 
In case (C.1), $p$ denotes a prime.\\
(ii) For all other equivalence classes, $G=\bbz_n$.
\end{theorem}

\noindent{\bf Comments}:\\
(1) Note that if none of $a, b$ or $c$ is prime to $n$ (but still assuming
$[n,a,b,c]=1$) then the curve
$y^n = x^a(x-1)^b(x+1)^c$ has full automorphism group $\bbz_n$.  It
is also interesting to note that the only cyclic cover with three
branch points that is Hurwitz (i.e. with full group of
automorphisms of order $84(g-1)$, where $g$ is the genus) is the Klein
curve of genus $3$ (case C.2). \\
(2) A discussion of the group structures and presentations as well as 
explicit descriptions of some of the various $G$ actions on $C$ is
given in \S\S\ 4, 5. In fact a nice aspect of this work is to
give explicit equations for some of the actions.\\
(3) The surfaces A.1 and A.2 are hyperelliptic, with the hyperelliptic
involution generating the
$\bbz_2$ factor in their automorphism group. Both surfaces are of ``Fermat
type'' and are special cases of {\bf Theorem \ref{Fermat}} below.  The surface
$y^{2m}=x(x-1)(x+1)^{2m-2}$ in A.2 (here $n=2m$) is isomorphic
to $y^{2m}=x^2-1$ and is known as the Accola-Maclachlan surface (see
{\bf Example \ref{accola-maclachlan}}). \\
(4) The curve $C: y^{2m}=x(x-1)^{m-2}(x+1)^{m+1}$ in B.2 ($n=2m$,
$4|m$) is on the other hand known as the Kulkarni surface (with
symmetries given in {\bf Example \ref{kulkarni}}). Note that the surface in
B.3 is a Kulkarni surface as well (corresponding to $m=4$) and has an
extra $\bbz_3$ action in its automorphism group. The curves B.1-B.3
are related by the fact that $a=1$, and either 
$b^2\equiv 1\ (mod\ n)$
or $c^2\equiv 1\ (mod\ n)$.\\
(5) Some of the symmetry groups above are given by central extensions of
cyclic groups by polyhedral groups. This is true for cases D.1-E.3 for
instance. $GL(2,3)$, the group of non-singular $2\times 2$ matrices over
$\bbf_3$, is isomorphic to  $(central\ \bbz_2):S_4$.  The geometry behind this is illustrated in part in {\bf Theorem
\ref{central}} below (with further details and examples given in \cite{KSS}).\\
(6) The last three surfaces E.1-E.3 are related as follows.
Taking the quotient of the surface $y^{24}=x(x-1)^4(x+1)^{19}$ (case E.3) by
$\bbz_{3}\subset$ (central $\bbz_{6}$), 
we get the surface \mbox{$y^8=x(x-1)^4(x+1)^3$}
(since $19\equiv 3\ (mod\ 8)$) which is the surface in E.1. Its group of
automorphisms is then $\bbz_2=\bbz_6/\bbz_3$ extended by $S_4$ and
this is same as $GL(2,3)$ . Taking the quotient of the surface E.3 by 
$\bbz_2\subset$ (central $\bbz_6$)
we get E.2.\\
(7) The surface E.1 turns out to be {\sl unique} in the sense that
there is a unique curve of genus 2 which affords an action by $\bbz_8$
(the cyclic action and the genus completely determine the surface;
see \cite{Kl}). This also happens for the surface $y^5=x(x-1)$ 
(see {\bf Remark 3}).

A surface is called {\sl belyi} if it branches over three points on
the line. Such curves, by a remarkable theorem of Belyi, have the
property of being isomorphic to curves defined over $\bar\bbq$
(see \cite{CIW}). A special class of belyi covers are the {\sl Lefschetz}
surfaces (see \cite{RR}) with equation
$$y^p = x^a(x-1),~0<a< p-1,~~p~\hbox{prime.}$$ 
These were originally studied in \cite{L}. In fact we can assume that
$1\leq a< \frac{p-1}{2}$ (i.e. if $a>\frac{p-1}{ 2}$, replace $a$ by
$p-a-1$, and if $a=\frac{p-1}{2}$, replace $a$ by $1$. This doesn't
change the isomorphism type of the curve; see \S5).  

As a corollary of the classification in {\bf Theorem \ref{maintheorem}} we
are able to recover (in particular) the following calculation of
Lefschetz which originally used an impressive mix of rational
functions, abelian varieties and divisor theory.

\bth
Let $C: y^p = x^a(x+1)$ be a Lefschetz surface and let $G:=$Aut$(C)$. 
Then\\
(1) If $a=1$, $G$ is cyclic of order $2p.$\\
(2) If $p=7$ and $a=2$, then $G=PSL(2,7)$ is the simple group of order
$168$.\\
(3) If $p\equiv 1\ (mod\ 3)$, $p>7$ and $1+a+a^2\equiv 0\ (mod\ p)$, 
then $G$ is the unique non-abelian group of order $3p$.\\
(4) For all other cases, $\aut (C)=\bbz_p$.
\end{theorem}

It is possible in special cases to solve for the automorphism group of
cyclic covers that are not belyi (no systematic technique is
known). The situation is particularly interesting for the Fermat
curves $F: y^n + x^d = 1$, with $d\leq n$. For $d\geq 4$, the $\bbz_n$
action in this case is not anymore uniformized by a triangle group and
hence our methods do not apply directly. There is however a way to get
around this by not constraining ourselves to $\bbz_n$ but by
uniformizing the entire $\bbz_d\oplus\bbz_n$ action on $F$. When this
is done, triangle groups appear again and the classification of
Conder, Bujalance and Cirre \cite{BCC} can be applied to the
situation. We summarize our calculations in this case.

\bth \label{Fermat}
Let $F$ be the surface given by $y^n+x^d=1$, where $4\le d \le n,$
and let Aut$(F)$ be its group of automorphisms.\\
(1) If $d=n$, then Aut$(F)= (\bbz_n\oplus\bbz_n)\ltimes S_3$ 
where $S_3$ is the symmetric group on 3 letters.  \\ 
(2) If $d$ does not divide $n$ then Aut$(F)=\bbz_d\oplus\bbz_n$. \\
(3) If $d|n$, $d<n$, Aut$(F)$ is the central $\bbz_d$ extension by the dihedral
group $D_{2n}$, given by the presentation
$$ \hbox{Aut} (F) = 
\big < s,t,u\ \bigm |\ s^d=t^n=u^2= [s,t]=[s,u]=1, (ut)^2=s^{-1}\big >$$
\end{theorem}

The special cases $d=2,3$ are covered in {\bf Examples 
\ref{special7}, \ref{special8}}. The automorphism group of the classic 
Fermat curve $x^n+y^n=1$ is of
course well-known (see \cite{Sh} or \cite{T}). The  involution $u$ 
in the case $d|n$
is fairly easy to describe. If $n=md$ then
$u:(x,y)\mapsto \left(\frac{x}{y^m}, \frac{1}{y}\right)$ is an
involution acting on the surface $y^n+x^d+1=0$, which is isomorphic to $F$.

The ideas we use in the proofs of these theorems are classical.  The
first basic idea in determining $\aut (C)$ is to associate to a Galois
cover $\beta: C\fract{G}{\lrar}\bbp^1$ a short exact sequence of
groups. Write $C={\bf U}/\Pi$ where $\Pi$ is a torsion free (Fuchsian)
group acting fixed point freely on $\bf U$, the upper half plane (see
\S2).  Then there is a Fuchsian group $\Gamma\subset PSL_2(\bbr )$
and a short exact sequence
$$1\lrar \Pi\lrar \Gamma\fract{\theta}{\lrar} G\lrar 1\leqno{(E)}$$ 
where $\theta$ is an epimorphism with torsion free kernel (or {\bf
skep} for surface kernel epi). We say that the sequence uniformizes
the action. Let $N(\Pi )$ be the
normalizer of $\Pi$ in PSL$_2(\bbr )$.  Then $N(\Pi )$ is itself a
Fuchsian group and $\hbox{Aut}(C) = N(\Pi )/\Pi$. When $C$ is
a cyclic cover of the line with group $\bbz_n$, then
$\Gamma_n:=\Gamma$ has a very special form. It has signature $\left(0\
|\ \frac{n}{[n,a_1]},\ldots, \frac{n}{ [n,a_k]}\right)$, see \S3.

The core of our work consists first of all in analyzing all possible
skeps $\theta_n: \Gamma_n\lrar\bbz_n$ and then seeking extensions to
larger Fuchsian groups $\theta'_n: \Gamma'_n\lrar G'$, with
$\bbz_n\subset G'$ and ker($\theta_n'$)=ker($\theta_n )$ = $\Pi$.  If
$\Gamma'$ happens to be finitely maximal, then necessarily $G'=\aut
(C)$. This method of extendability through skeps works very
well for {\sl triangle groups} $\Delta$ with signature $(0\ |\
m_1,m_2,m_3 )$, and hence the content of the main theorem.

Finally, we observe that a type of stability occurs:
when the number of branch points is large enough, the cyclic
$\bbz_n$ action necessarily normalizes.

\bth\label{central} 
Let $C:y^p= f(x)$, where $p$ is a prime and $f(x)$ is a polynomial
with $r$ distinct roots, $r> 2p$. Then 
the automorphism group of $C$ is an
extension of $\bbz_p$ by a polyhedral group.
\end{theorem}

A good example is already illustrated by the Fermat curve
$y^n+x^d=1$, $d|n$, with automorphism group $\bbz_d: D_{2n}$.

{\sc Note Added}: While this paper was under revision, we were told
that some recent calculations of Bujalance, Cirre and Turbek
overlapped with ours.


\vskip 10pt
{\bf\Large\S2 Group Actions and Extensions of Fuchsian Groups}

Let $C$ be a Riemann surface of genus $g\geq 2$. By the uniformization
theorem, $C$ is the quotient of a discrete (torsion free) group $\Pi$
acting fixed-point freely on the upper half plane.  The holomorphic
structure on $C$ is induced by the quotient map and so depends not
only on the abstract isomorphism class of the group $\Pi$ but also on 
the way this group embeds in PSL$_2(\bbr )$.

By a Fuchsian group we mean any discrete subgroup of PSL$_2(\bbr
)=$Isom $(\bf U )$, where $\bf U$ is the upper-half plane. 
A Fuchsian group with compact quotient has presentation (see \cite{S}):
$$\Gamma = \left< x_1,\ldots, x_r; a_1, b_1, \ldots , a_g, b_g~|~
x_1^{m_1} = \cdots = x_r^{m_r} = x_1\cdots x_r [a_1,b_1]\cdots
[a_g,b_g] = 1\right>$$ 
 Here $g$ is the genus of ${\bf U}/\Gamma$. We encode this
in a ``signature'' $\sigma (\Gamma ) = (g~|~m_1,\ldots, m_r )$ 
and when $g=0$ we write $\sigma (\Gamma ) = (m_1,\ldots, m_r)$.
In this case the action of the elliptic element $x_i\in\Gamma$ 
on $\bf U$ is given by rotation by ${2\pi/ m_i}$ about a fixed point.  

When $r=0$,  the group $\Gamma$ is torsion free and corresponds
to the universal covering group of a Riemann surface. The signature in
this case is $\sigma (\Gamma )= (g~|~-)$ and $\Gamma$ is isomorphic to
the fundamental group of the genus $g$ surface ${\bf U}/\Gamma$.  The
group $\Gamma$ is of genus zero if $g=0$. 

Let $\Gamma$ be the uniformizing Fuchsian group for the action of $G$
on $C$, as in the introduction,
and let $\Gamma^{\prime}$ be another Fuchsian
group containing $\Gamma$ with finite index. We say that the action of
$G$ extends to $G^{\prime}$ if there is a commuting diagram
$$\begin{array}{ccccccccc}
1&\lrar&\Pi&\lrar &\Gamma&\fract{\theta}{\lrar}&G&\lrar& 1\\
&&\decdnar{=}&&\decdnar{\nu}&&\decdnar{\mu}\\
1&\lrar&\Pi&\lrar &\Gamma^{\prime}&\fract{\theta}{\lrar}&G^{\prime}&\lrar& 1
\end{array}
$$
where $C\cong {\bf U}/\Pi$ and $\mu,\nu$ are inclusions. Now clearly 
$G\subset G'\subset $Aut$ (C)$.
If $\Gamma'$ is maximal then it must coincide with the normalizer
of $\Pi$ in PSL$_2(\bbr )$ and hence

\ble
If $\Gamma^{\prime}$ is maximal then $\hbox{Aut}(C) = G^{\prime}.$
\ele

From this viewpoint, finite extendability of Fuchsian groups is a
necessary step in computing automorphism groups of Riemann surfaces.
Given an arbitrary Fuchsian group $\Gamma$ it is usually not possible to find a proper
inclusion  $\Gamma\subset\Gamma'$ of finite index, in which case  
$\hbox{Aut}(C) = G.$
The geometry of the fundamental domain of $\Gamma$
plays a seminal role in the existence or non-existence of such
extensions (with an exception for the triangle groups as is explained
below).  In general, some embeddings of the abstract group $\Gamma$ as a Fuchsian group
will admit extensions and others will not.

It turns out however that for a certain class of Fuchsian groups
$\Gamma$, every monomorphism $\rho:\Gamma\lrar$ PSL$_2(\bbr )$ extends
to a monomorphism $\rho^{\prime}: \Gamma^{\prime}\lrar$ PSL$_2(\bbr )$ 
for some
$\Gamma^{\prime}$ of finite index. Such groups are said to have ``non-finitely
maximal signature'' and a list of them is given below (we only look at
genus 0 curves).

\bth (Greenberg, Singermann):\label{table}
The only genus zero Fuchsian groups $\Gamma$ with non-finitely maximal signature
are those on the following list
$$\begin{tabular}{|c|l|c|c|c|}\hline
&$\hspace{0.8in}\Gamma$&$\Gamma'$&$|\Gamma':\Gamma|$\\
\hline
1& $(n, n, n), n\geq 4$&$(3, 3, n )$&$3$\\
2& $(n, n, n), n\geq 4$ &$(2, 3, 2n)$&$6$\\
3& $(n, n, m), n\geq 3, n+m\geq 7$&$(2, n, 2m)$&$2$\\
A& $(n,n,n,n), n\geq 3$&$(2,2,2,n)$&$4$\\
B& $(n,n,m,m), n+m\geq 5$&$(2,2,n,m)$&$2$\\
\hline
4& $(7,7,7)$&$(2,3,7)$&$24$\\
5& $(2,7,7)$&$(2,3,7)$&$9$\\
6& $(3,3,7)$&$(2,3,7)$&$8$\\
7& $(4,8,8)$&$(2,3,8)$&$12$\\
8& $(3,8,8)$&$(2,3,8)$&$10$\\
9& $(9,9,9)$&$(2,3,9)$&$12$\\
10& $(4,4,5)$&$(2,4,5)$&$6$\\
11& $(n,4n,4n), n\geq 2$&$(2,3,4n)$&$6$\\
12& $(n,2n,2n), n\geq 3$&$(2,4,2n)$&$4$\\
13& $(3,n,3n), n\geq 3$&$(2,3,3n)$&$4$\\
14& $(2,n,2n), n\geq 4$&$(2,3,2n)$&$3$\\
\hline
\end{tabular}$$
The first column in this table shows $\Gamma$ (or rather its
signature), the second an extension $\Gamma^{\prime}$ of $\Gamma$ and the
third the index of the extension. Only the first five
extensions are normal. 
\end{theorem}

We refer to the table above as the GS table. A signature (with
abstract group $\Gamma$) is therefore finitely maximal if for some
embedding of $\Gamma$ in PSL$_2(\bbr )$, the Fuchsian group so
obtained is finitely maximal.

Now triangle groups $\Delta$ with signature $(m_1, m_2, m_3 )$ have
the following special properties:\\ 
(1) All embeddings of $\Delta (m_1, m_2, m_3 )$ in PSL$_2(\bbr
)$ are conjugate.\\
(2) If $\Delta\subset \Gamma$ is a finite index extension of a
triangle subgroup, then $\Gamma$ is itself a triangle group (\cite{B},
theorem 10.6.5).  \\
The first property means that the existence of a finite index
extension of $\Delta (m_1, m_2, m_3 )$ in PSL$_2(\bbr )$ does not 
depend on the way the group embeds. So a triangle group 
either always extends or never extends.

\begin{figure}[htb]
\begin{center}
\epsfig{file=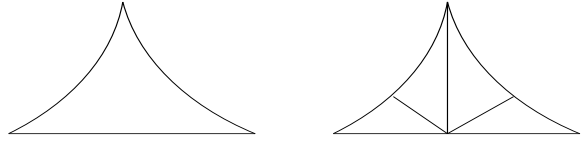,height=0.7in,width=3in,angle=0.0}
\caption{Some inclusions $\Delta(n,2n,2n)\subset\Delta(2,4,2n)$}\label{incl2}
\end{center}
\end{figure}

{\bf The embeddings}: There is nice geometry behind 
the embedddings described in the GS table. Often it is
possible to deduce inclusions of triangle groups $\Gamma\subset\Gamma'$
by subdividing the triangle associated to $\Gamma$ into $n$ copies of the 
triangle associated to $\Gamma'$, where $n$ is the index of $\Gamma$ in 
$\Gamma'$. We
illustrate this in the case $\Delta (n, 2n, 2n)\hookrightarrow\Delta
(2,4, 2n)$ (see Figure ~\ref{incl2}).

The left hand diagram is a triangle with angles $\pi/2n$, $\pi/2n$
(bottom) and $\pi/n$ (top). We can bisect it by 
dropping a perpendicular from the top vertex.  This gives us $2$ 
copies of a triangle with angles $\pi/2,\pi/2n,\pi/2n$  and
illustrates an index 2 inclusion (there are $2$ of them):
$\Delta(n,2n,2n)\subset\Delta(2,2n,2n).$ 
Constructing 2 more perpendiculars
gives $4$ copies of the triangle with angles $\pi/2,\pi/4,\pi/2n$
and $4$ inclusions 
 $\Delta (n, 2n, 2n)\hookrightarrow\Delta
(2,4, 2n).$
In fact, by an appropriate choice of elliptic generators
$x_1$, $x_2$ and $x_3$ for $\Delta (n, 2n, 2n)$ and 
$y_1$, $y_2$ and $y_3$ for $\Delta (2,4,2n)$ an embedding 
$\Delta (n, 2n, 2n)\hookrightarrow\Delta
(2,4,2n)$ can be described by
$$x_1\mapsto y_3^{-1}y_1y_3^2y_1y_3,~~x_2\mapsto
y_1y_2^{-1}y_3y_2y_1,~~x_3\mapsto y_3$$

The following assertion (whose proof we skip) relies on a tedious
case-by-case study.

\bpr
Any two embeddings $\Gamma\hookrightarrow\Gamma^{\prime}$ in the GS table
differ by an automorphism of $\Gamma^{\prime}$.
\epr

{\bf Corollary}: If $\theta: \Gamma\lrar G$ extends to $\Gamma^{\prime}\lrar
G^{\prime}$ for some embedding $\Gamma\subset\Gamma^{\prime}$, $\Gamma$ 
as in the GS
table, then $\theta$ also extends for all other embeddings.


\vskip 10pt
{\bf\Large\S3. Cyclic Galois Coverings of the Line}

Consider the affine curve 
$$C: y^n = p(x),\  \mbox{where}\ p(x)\ \mbox{is a polynomial of
degree}\ m$$ 
Write $f(x,y) = y^n - p(x)$.  The points in $\bbc^2$ where the curve
has singularities (i.e. where $\frac{\partial f}{\partial x}= 0 =
\frac{\partial f}{\partial y}$) occur at points $(x,y)=(a,0)$ where $a$
is a multiple root of $p(x)$. To get a compact Riemann surface from
the affine model of $C$, we need to remove these singularities (by
``adding'' ramification points), and then compactify in
$\bbc^2\subset\bbp^2$ (that is adding points at infinity in a
prescribed way). The points at infinity may themselves be singular so
this needs to be analyzed as well.

Consider the projection $\pi: C\lrar\bbc$, which sends $(x,y)\mapsto
x$. Obviously $\pi$ is going to be regular everywhere but at the roots
of $p$ (these are the branch points and their preimages we call the
ramification points). Let $a$ be a root of $p$ and write $y^n = (x-a
)^kq(x)$. If \mbox{$k\equiv j\ (mod\ n)$,} then by the birational
transformation $(x,y)\mapsto (x, y/ (x-a)^l)$, where \mbox{$k=ln+j$,} we see
that $C$ is biholomorphic to the curve $y^n = (x-a)^jq(x)$. So we
might as well assume to begin with that $C: y^n = (x-a)^kq(x)$ with
$k<n$.

In a sufficiently small neighborhood of $a$, $q(x)$ is non-zero and by
picking an $n$-th root branch we can absorb it with the term in $y$ so
that the curve (around the singularity $x=a$) has the equation $w^n =
(x-a)^k$ (where $w^n=y^n/q(x)$).  Let $l= [n,k]$, the greatest common
divisor, and let $n=ln', k=lk'$.  Then
$$w^n - (x-a)^k= \prod_{i=0}^{l-1}(w^{n'} - \zeta^i(x-a)^{k'} )$$
where $\zeta$ is a primitive $l$-th root of unity.  Each term $w^{n'} -
\zeta^i(x-a)^{k'} = 0$ has a {\sl removable} singularity at 
$(a,0)$ because $[n',k']=1$; see \cite{Mi}, p. 71.  So by adding
a point there, locally the sheets of $C': w^{n'} - \zeta^i(x-a)^{k'} =
0$ come together in a smooth way. Going back to $C$, adding a
(ramification) point for each of the terms 
$w^{n'} -\zeta^i(x-a)^{k'} $ gives us a smooth surface. In total we have added
$[n,k]$ such points.

The same technique, after a slight change of variables, allows one to
study what happens in the neighborhood of infinity (see \cite{Mi}, p. 73).
Write $C: y^n = c(x-a_1)^{k_1}\ldots (x-a_m )^{k_m}$. Then it can be
checked that if $\sum k_i\equiv 0\ (mod\ n)$, one needs to add $n$ points
$x_1,\ldots, x_n$ to $C$ over $\infty$ in order that the projection $(x,y)\mapsto x$ can
be extended to a ramified covering (also denoted by $\pi$) from the
compactified $C$ to $\bbp^1 = \bbc\cup\infty$. There is no branching 
over $\infty$ in this case.  On
the other hand if $\sum k_i = qn + r$, $0<r<n $, then one needs to add
$[n,r]$ points over infinity to compactify the surface (in which case
$\pi$ is branched over $\infty$).

The branch point data for $\pi: C\lrar\bbp^1$ is now totally
explicit. At each ramification point over the branch point
$(y=0, x=a_i )$, there are $n/[n,k_i ]$ sheets coming
together (this is the multiplicity of the ramification point).
If we assume $\sum k_i\equiv 0\ (mod\ n)$ (to avoid branching over
$\infty$), then we can immediately read off the genus of $C$ from 
the Riemann-Hurwitz formula
$$2-2g = n\left( 2 - \sum_{i=1}^m\left( 1- \frac{[n,k_i]}{ n}\right)\right )$$
$$\Longrightarrow \ g = 
\frac{1}{2}\left( 2 +n(m-2)- \sum^m_{i=1} [n,k_i]\right)$$
Note that if all the $k_i$ are prime to $n$ then 
$g = \frac{(n-1)(m-2)}{2}$ if $\sum k_i\equiv 0\ (mod\ n)$.

\bpr \label{cyclicovers}
Consider the compact Riemann surface $C\subset \bbp^2$ associated
to the affine curve $C: y^n = (x-a_1)^{k_1}\ldots (x-a_m
)^{k_m}$, $m\geq 3$, and assume $\sum k_i\equiv 0\ (mod\ n)$. The curve is
irreducible if $[n,k_1,\ldots, k_m]=1$. Moreover, its
genus is given by the formula
$g = \frac{1}{ 2}\left[ 2 + (m-2)n -\sum [n,k_i]\right]$.
The cyclic action of $\bbz_n$ on the curve is
uniformized by an exact sequence
$$1\lrar\Pi\lrar \Gamma \left( \frac{n}{[n,k_1]},
\ldots, \frac{n}{[n,k_m]} \right)\fract{\theta}{\ra 3}\bbz_n\lrar 1$$
where $\theta$ is defined on the elliptic generators $x_i\in \Gamma$
(of order $n/[n,k_i]$) by $\theta (x_i) = T^{k_i}$, $T$ being 
a generator of $\bbz_n$.
\epr

{\sc Proof}:
It follows from a theorem of Capelli and Kneser (\cite{Sz}, p. 92) 
that the binomial $y^n-f(x)$ is irreducible over $\bbc (x,y)$, and hence
over $\bbc [x,y]$, if $[n,k_1,\ldots, k_m]=1$. 

Now ${\bf U}\lrar {\bf U}/\Gamma = \bbp^1$ is ramified exactly at
the branch points of $\pi:C\lrar\bbp^1$ (given by the $a_i$). The
periods of $\Gamma$ coincide with the order of the stabilizers at the
ramification points. At each of these points (lying over $a_i$ say),
there are $n/[n,k_i]$ sheets coming together and the stabilizer
group (at any of these points) is necessarily cyclic (of that
order). This yields the assertion about the form of $\Gamma$.

Finally, to analyze $\theta$ we need to understand how the sheets come
together to give the smooth surface $C$ (this is encoded in the
monodromy representation of the associated unbranched or \'{e}tale
cover). Let $a$ denote one of the branch points and $k$ the corresponding
exponent.  Then  $C$ is locally given
by $y^n = (x-a)^kg(x)$, where $g(a)\ne 0.$ Choose a basepoint $x_0\in\bbc-\{a_1,\ldots,
a_m\}$ and let $w_a$ denote a smooth simple closed curve based at
$x_0$ and going once around $a$ (chosen so that the other branch
points are in its exterior). At the point $x_0$ we can choose an
arbitrary branch $y_1$ of the multi-valued function $y=f(x)$ and then
choose the other branches by $y_j=\zeta^{j-1}y_1$, $1\leq j\leq n$,
where $\zeta = e^\frac{2\pi i}{ n}$.

Analytic continuation of any germ $y$ at $x_0$ once around $w_a$ takes
us to $\zeta^ky$, and therefore the cycle decomposition of the
monodromy representation at $a$ is given by
$$\pi_a = (y_{i_1},\zeta^ky_{i_1},\zeta^{2k}y_{i_1},\ldots )
\times\cdots\times (y_{i_{[n,k]}},\zeta^ky_{i_{[n,k]}},
\zeta^{2k}y_{i_{[n,k]}},\ldots )$$
(there are $[n,k]$ such factors). Here the $y_{i_r}$ are
representatives of the cycles. let $\tilde C$ denote the unbranched
Riemann surface associated to $y^n=\prod^m(x-a_i)^{k_i}$. The
unbranched covering ${\tilde C}\lrar \bbc-\{a_1,\ldots, a_m\}$ 
yields a short exact sequence
$$1\lrar\pi_1({\tilde C},{\tilde x}_0)\lrar
\pi_1(\bbc-\{a_1,\ldots, a_m\}, x_0)\fract{\phi}{\lrar}\bbz_n\lrar 1$$
Now $\pi_1(\bbc-\{a_1,\ldots, a_m\}) = <w_1,\ldots, w_m\ |\ 
w_1\cdots w_m = 1>$. The action of $\phi (w_i )$ on $\tilde C$ is determined
by what happens over $x_0$ and this is determined by the monodromy 
$\pi_{a_i}$. Since $\pi_{a_i}(y)=\zeta^{k_i}y$ for any germ $y$ based
at $x_0$, we have $\phi (w_i)=T^{k_i}$.
\hfill\za

{\bf Remark 1}: In \cite{H}, Harvey has given necessary and sufficient
conditions for a skep $\theta:\Gamma (m_1,\ldots, m_r)\lrar\bbz_n$ to
exist. His main condition is that
$$n = \hbox{lcm}(m_1,m_2, \ldots, m_r )=
\hbox{lcm}(m_1,\ldots, \hat m_i, \ldots, m_r ),\ \forall i$$
where $\hat m_i$ means deleting the $i$-th entry. The existence of the
uniformizing sequence in {\bf Proposition \ref{cyclicovers}} 
 shows that if $k_i < n$, $\sum
k_i\equiv 0\ (mod\ n)$ and $[n,k_1,\ldots, k_m]=1$, then necessarily

$$ n = \hbox{lcm}\left(\frac{n}{ [n, k_1]},\ldots, \frac{n}{ [n, k_m]}\right) 
= \hbox{lcm}
\left(
\frac{n}{[n, k_1]},\ldots, \widehat{\frac{n}{[n, k_i]}}, 
\ldots, \frac{n}{[n, k_m] }
\right)
$$
for all $1\leq i\leq r $. This can be checked directly of course.

\ble \label{model} Let $C: y^n = \prod^m_{i=1} (x-a_i)^{k_i}$ and let
$l$ be prime to $n$, $k^{\prime}_i\equiv lk\ (mod\ n)$. Then $C$ is birationally
equivalent to the curve $C': y^n = \prod^m_{i=1} (x-a_i)^{k'_i}$.
\ele

{\sc Proof}: Consider the curve $C_l: y^n =\prod^m_{i=1} (x-a_i)^{lk_i}$ and
let $\psi: C\lrar C_l, (x,y)\mapsto (x,y^l )$.  We check that this map
is 1-1 and hence a biholomorphism. Suppose $(x_1,y_1^l)=(x_2,y_2^l)\in
C_l.$  Then $x_1=x_2$ and $y_1^l=y_2^l$. But $(x_1,y_1),(x_2,y_2)\in
C$ then implies that $y_1^n=y_2^n$.  Now $[l,n]=1$ and hence there is $k<n$
such that $kl\equiv 1\ (mod\ n)$. Write $kl=\beta n+1$, then
$$y_1^{lk}=y_2^{lk}\Longrightarrow y_1^{\beta n}y_1 = y_2^{\beta
n}y_2\Longrightarrow y_1=y_2$$ 
This shows that $\psi : C\cong C_l$. Finally, one
absorbs all $n$-th powers into the $y$ term to show that
$C'\cong C_l$.\hfill\za

{\sc Automorphism groups}: The Galois group in this case is $G=\bbz_n$
and the uniformizing sequence is given by
$$1\lrar \Pi\lrar \Gamma\fract{\theta}{\lrar} \bbz_n\lrar 1$$ As
explained in \S2, we seek extensions to $\Gamma'\lrar G'$ for some
finite group $G'$ containing $\bbz_n$. We restrict our attention to those
$\Gamma,\Gamma^{\prime}$ in the GS table. 

Now a case-by-case study of all such possible extensions (exhausting
the GS table) has already been carried out by Bujalance and
Conder. Denote by $x_1, \ldots, x_m$ the generators of $\Gamma
(n_1,\ldots, n_m)$ and let $z_1=\theta (x_1),\ldots, z_m=\theta (x_m)$
be their images in $\bbz_n$. Also write $T$ for a generator of
$\bbz_n$.  Their main result is:

\bth \label{conder-bujalance}
\cite{BC} Suppose $\bbz_n$ acts on $C={\bf U}/\Pi$.
Then, in the following cases
the action of $\bbz_n$ can be extended to an action of a larger
group and so Aut$(C)\neq \bbz_n$:
\begin{enumerate}
\item sig$(\Gamma ) = (n,n,n,n )$, $z_1 =T$, 
$z_2 =T^a, z_3=T^b$, $z_4 =T^c$ where $abc\equiv 1\ (mod\ n)$,
$a^2\equiv b^2\equiv c^2\equiv 1\ (mod\ n)$ and
$1+a+b+c\equiv 0\ (mod\ n)$.
\item sig$(\Gamma ) = (n,n,m, m)$, $n+m\ge 5$, in which case $\bbz_n$
extends to  the dihedral group 
$D_{2n}= \left< u,v\ |\ u^2=v^n =(uv)^2=1\right>$
 of order $2n.$
\item sig$(\Gamma ) = (n,n,n )$ if $n\geq 4$ and
$\bbz_n$ has an  automorphism of order $3$
permuting $z_1, z_2, z_3$.
\item  sig$(\Gamma ) = (n, n, m)$ where $n\geq 3$, $m|n$,
and $n+m\geq 7$, and either 
$z_1=z_2$ or there is a
transposition exchanging $z_1$ and $z_2$.
\item sig$(\Gamma ) = (3,4, 12)$, with
$\{z_1,z_2,z_3\}=\{T, T^{3}, T^{-4}\}$. 
\end{enumerate}
\end{theorem}

The converse of this theorem is true if we restrict attention to triangle 
groups.  That is Aut$(C)=\bbz_n$ if the action is uniformized by a triangle
group not covered by cases (3), (4) or (5).

The discussion in \cite{BC} is actually sufficient to give the automorphism
group in each situation. Starting with this classification
we now seek extensions of
$$\theta: \Gamma \left(\frac{n}{ [n,k_1]},\ldots, \frac{n}{ [n,k_m]}\right)
\lrar\bbz_n,\ m=3,4$$
where $1\leq k_i<n$, $[n,k_1,\ldots, k_m]=1$ and 
$k_1+\cdots +k_m\equiv 0\ (mod\ n)$.  As an example, the next corollary 
follows immediately from {\bf Proposition
\ref{cyclicovers}} and {\bf Theorem \ref{conder-bujalance}}.

\bco Let $C: y^n = (x-a_1)^{k_1}\cdots (x-a_4)^{k_4}$ with the above 
assumptions.  Assume in addition that
 $[n,k_1]=[n,k_2]$ and $[n,k_3]=[n, k_4]$.  Then the
dihedral group $D_{2n}$ acts on $C$.  \eco

Finally we record the following observation ({\bf Theorem \ref{central}} in
the introduction ). 

\bth
Let $C:y^p= f(x)$, $p$ prime and $f(x)$ a polynomial with $r$ distinct
roots, $r> 2p$. Then the automorphism group of $C$ is an extension of
$\bbz_p$ by a polyhedral group.
\end{theorem}

\noindent{\sc Proof}: Let  $\nu (\tau )$ denote the number of fixed
points of an  automorphism $\tau: C\lrar C$.  The
following  result (see \cite{FK}, p. 245]) is crucial.  Let $g:
C\lrar\bbp^1$ be a meromorphic map and let $\tau: C\lrar C$ be an
automorphism with $\nu (\tau)> 2\deg g$. Then necessarily $g$ is
invariant under $\tau$ (i.e. $g\circ \tau =g$).

Suppose  $\bbz_p=<T>$ is the cyclic group
of automorphisms acting on $C$ by $T(x,y)=(x,\zeta y )$.  Apply the above
result to the quotient map $g:C\lrar C/<T>\cong \bbp^1$ and 
$\tau=RTR^{-1}$, $R\in\hbox{Aut}(C)$. The degree of $g$ is $p$ and therefore
$RTR^{-1}=T^k$ for some $k.$  Hence $<T>$ is normal in $\hbox{Aut}(C)$
and the result follows by standard covering space theory. \hfill\za


\vskip 10pt
{\bf\Large\S4 Belyi Curves and Lefschetz Surfaces}

By a ``generalized Lefschetz'' curve we mean an algebraic curve
given by the affine equation
$$L: y^n = x^a(x-1)^b(x+1)^c,\ 1\le a,b,c\le n-1, \ 
a+b+c\equiv 0\ (mod\ n)$$ 
(the classic Lefschetz case corresponds to $n=p$, a prime). It has genus
$$g = \frac{1}{2}(2+n-[n,a]-[n,b]-[n,c])$$
according to {\bf Proposition \ref{cyclicovers}}. Note that when
$a,b,c$ are prime to $n$, the
genus is $g=\frac{1}{2}(n-1)$.
As was discussed in the introduction, generalized Lefschetz curves are
belyi with belyi map $\beta: L\to\bbp^1, (x,y)\to x$, which is
branched over $0,1$ and $-1$.  The $\bbz_n$ action is uniformized by
a skep
$$\theta: \Delta\left(\frac{n}{[n,a]},\frac{n}{ [n,b]},\frac{n}{
[n,c]}\right)\lrar\bbz_n$$ 
The main purpose of this section and the next is to analyze the
extendability of such skeps.  This means looking closely into cases
(3), (4) and (5) of {\bf Theorem \ref{conder-bujalance}}. We write
$z_1, z_2, z_3$ for the images of the elliptic generators 
$x_1,x_2,x_3$ under
$\theta :\Delta (m_1,m_2,m_3)\lrar\bbz_{n}$.

{\underline{Case (3)}}: This is the case 
$\Delta \left({n/ [n,a]},{n/[n,b]},{n/ [n,c]}\right) =\Delta (n,n,n)$, 
that is all of $a,b,c$ are prime
to $n$. By {\bf Lemma \ref{model}} we can choose $a=1$. The only
possible extension in this case is the degree 3 extension of $\Delta
(n,n,n)\lrar\bbz_n$ to $\Delta (3,3,n)\lrar G'$ and this occurs if there
exists an automorphism of $\bbz_n$ of order 3 cyclically  permuting
$z_1$, $z_2$ and $z_3$ (\cite{BC}, case N6).  Let $\tau$ be such an
automorphism and write $\tau (T)=T^k$ for some $k$ prime
to $n$. Then $k^3\equiv 1\ (mod\ n)$. On the other hand we know that $z_1=T,
z_2= T^{b}$ and $z_3=T^{c}$.  Since $\tau (z_i)=z_{i+1}$ (index modulo
3) we find that
$$a=1,\ b=k,\ c\equiv k^2\ (mod\ n), \ [k,n]=[c,n]=1$$
and so necessarily $1+k+k^2\equiv 0\ (mod\ n)$ (and $n$ is odd).

Note that if $n=p_1^{n_1}\ldots p_k^{n_k}$ is the prime decomposition
then $\aut (\bbz_n )\cong \prod_{i=1}^k\aut (\bbz_{p_1^{n_i}})$. On
the other hand $\aut(\bbz_{p^{n}})=\bbz_{(p-1)p^{n-1}}$ if
$p$ is odd.  So if $3| |$Aut$\bbz_n |$, then either $9|n$ or there is an
odd prime $p$ such that $p|n$ and $p\equiv 1\ (mod\ 3)$.  Since we must
have $1+k+k^2\equiv 0\ (mod\ n)$, $9|n$ means $9|(1+k+k^2)$ which is not
possible for any integer $k$. Therefore, when there exists $p|n$, 
$p\equiv 1\ (mod\ 3),$ an index $3$ extension becomes possible
$$
\begin{array}{cccccc}
\Delta (n,n,n )&\lrar&\bbz_n\\
 \downarrow&&\downarrow\\
\Delta^{\prime}(3,3,n)&\lrar&\bbz_n\ltimes\bbz_3
\end{array}
$$
where $\bbz_n\ltimes\bbz_3=G'$ is a {\sl metacyclic} group of order
$3n$. In {\bf Example \ref{periodthree}} below we
explicitly construct the desired period 3 automorphism acting on the
curve $C: y^n=x(x-1)^k(x+1)^{k^2},$ when $n$ is odd
and $1+k+k^2\equiv 0\ (mod\ n)$.  This leads to cases $C_1$ and $C_2$ of
{\bf Theorem \ref{maintheorem}}.  It turns out that the skep 
$\Delta^{\prime}(3,3,n)\lrar\bbz_n\ltimes\bbz_3$ can be extended further
only in case $C_2.$ See the next section.

Case (4) being the longest to deal with, we first settle case (5).

{\underline{Case (5)}}: This is the case $\Delta
\left(n/[n,a],n/[n,b],n/[n,c]\right) =
\Delta(3,4,12)$, where $n=12.$  There is one possible extension (\cite{BC}, case T10)
of index 4
$$
\begin{array}{ccccccc}
\Delta(3,4,12)&\lrar&\bbz_{12}\\
 \downarrow&&\downarrow\\
\Delta^{\prime}(2,3,12)&\lrar&G^{\prime}
\end{array}
$$
with $G^{\prime}$ of order 48 given as a central extension 
$1\lrar \bbz_4\lrar
G^{\prime}\lrar A_4\lrar 1,$ where $A_4$ is the alternating group. Since
$[n,c]$ is prime to $12$ we can choose $c=1$ (by {\bf Lemma \ref{model}}).
According to (\cite{BC}, T10) the extension occurs if
$\{z_1,z_2,z_3\}=\{v, v^3, v^{-4}\},$ where $v$ is a generator of
$\bbz_{12}$.  Since $[12,a]=4$, $[12,b]=3$, it follows that $a=4$ or $8$
and $b=3$ or $9$. On the other hand $a+b+c\equiv 0\ (mod\ 12)$ and hence
(after permutation) $a=1, b=3, c=8$. Since $\Delta (2,3,12)$ is
maximal, we deduce that $G=$Aut$\left(y^{12}=x(x-1)^3(z+1)^8\right)$
has order $48$.  This leads to case $D_1$ of {\bf Theorem \ref{maintheorem}}.

{\underline{Case (4)}}: This is the case where $a$ and $b$ are prime to $n$, so
that $\sigma (\Gamma )$ is $(n,n, n/[n,c])$.  Note again, according to
{\bf Lemma \ref{model}} we can choose $a=1$.  Here $m =n/[n,c]$
and so from the GS table
there are the following subcases to consider: rows
2, 3, 4, 7, 9, 11 and 12.  According to {\bf Theorem \ref{conder-bujalance}},
cases 2 and 9 do not admit extensions, and since we have already considered
case 4, this leaves only cases 3, 7, 11 and 12.  In all cases, for an 
extension to occur the following condition is necessary: either 
$z_1=z_2$, or $z_1\neq z_2$ and there
is an involution of $\bbz_n$ exchanging $z_1$ and $z_2$.  See
{\bf Theorem \ref{conder-bujalance}}.

{\underline{Subcase 4.1}}: This refers to the index 2 extension
 in row 3 of the GS table
$$
\begin{array}{cccccc}
\Delta\left(n,n, n/[n,c]\right)&\lrar&\bbz_n\\
 \downarrow&&\downarrow\\
\Delta^{\prime}\left(2,n, 2n/[n,c]\right)&\lrar&G^{\prime}
\end{array}
$$
We distinguish two cases (see \cite{BC}, N8):\\
$\bullet$ Suppose $z_1=z_2$. Then $T^{a}=T^{b}$ and hence
$a=b=1$.  Thus $c=n-2$ and
the curve has the form
$C:y^n = x^{n-2}(x-1)(x+1),$ after permutation of $a,b,c.$   
$C$ admits the 
involution 
$(x,y)\mapsto (-x,-y)$ and so $G'=\bbz_2\oplus\bbz_n\subset\aut
(C)$.  If $n$ is odd then $[n,c]=[n,n-2]=1$ and thus
$\Delta\left(n,n, 2n/[n,c]\right)=\Delta\left(n,n, n\right)$.
The only possible extension is the composite
$$\Delta\left(n,n, n\right)\subset\Delta\left(2,n, 2n\right)\subset
\Delta\left(2,3,2n\right)$$
But this corresponds to row 2 of the GS table, and it is known that
no extension exists in this case.  This leads to case $A.1$ of
{\bf Theorem \ref{maintheorem}}. 

Now assume $n$ is even.  Then 
$$\Delta\left(n,n, n/[n,c]\right)=\Delta\left(n,n, n/2\right)\  \mbox{and}\ 
\Delta\left(2,n,2n/[n,c]\right)=\Delta\left(2,n,n\right)$$
But $\Delta\left(2,n,n\right)$ is not maximal.  
In fact it leads to the following 
possible extensions (from the GS table):
\begin{eqnarray*}\displaystyle
\Delta\left(n,n, n/2\right)&\subset&\Delta\left(2,n,n\right)\subset
\Delta\left(2,4,n\right)\\
\Delta\left(4,8,8\right)&\subset&\Delta\left(2,8,8\right)\subset
\Delta\left(2,3,8\right)\ \mbox{if $n=8$} 
\end{eqnarray*}
Each of these cases is considered below.\\
$\bullet$ Now assume $z_1\neq z_2$.  Then there is an involution $\tau$
interchanging $z_1$ and $z_2$, say $\tau
(T)=T^k$.  Since $\tau^2=1$, $\tau\neq 1$ we have $k^2\equiv 1\ (mod\ n)$, 
$k\not\equiv 1\ (mod\ n)$. We can take $k=b$, so the curve has the equation
$C: y^n = x(x-1)^k(x+1)^{n-k-1}$. The group $G'$ is a twisted
product of $\bbz_2$ with $\bbz_n.$  More relations can be deduced
between $n$ and $k$ (see {\bf Example \ref{twistedz2}}).  Again there are 
possible further extensions.  Let $m=n/[n,c].$
\begin{eqnarray*}\displaystyle
\Delta\left(n,n,n\right)&\subset&\Delta\left(2,n,2n\right)\subset
\Delta\left(2,3,2n\right)\ \mbox{if $n=m$}\\
\Delta\left(n,n,m\right)&\subset&\Delta\left(2,n,2m\right)\subset
\Delta\left(2,3,4m\right)\ \mbox{if $n=4m$}\\
\Delta\left(2m,2m,m\right)&\subset&\Delta\left(2,2m,2m\right)\subset
\Delta\left(2,4,2m\right)\ \mbox{if $n=2m$}\\
\Delta\left(4,8,8\right)&\subset&\Delta\left(2,8,8\right)\subset
\Delta\left(2,3,8\right)\ \mbox{if $n=8$}
\end{eqnarray*}
The first possibility corresponds to row 2 of the GS table, but this 
extension does not exist (see \cite{BC}). Each of the remaining cases is 
dealt with below.  Otherwise the extension is
maximal.  This leads to case $B.1$ of  {\bf Theorem \ref{maintheorem}}.

{\underline{Subcase 4.2}}: This refers to the index 12 extension
$$
\begin{array}{ccccc}
&\Delta (4,8,8)&\lrar&\bbz_8&\\
& \downarrow   & &\downarrow&\\
&\Delta (2,3,8)&\lrar&G'&
\end{array}
$$
In this case no further extension is possible if $z_1=z_2$, see 
(\cite{BC}, T4). Thus assume $z_1\neq z_2$.  Then we can choose
$a=1$, $b=5$ and $c=2$, see (\cite{BC}, T4), and an extension is
possible. The curve is given by 
$C: y^8= x(x-1)^2(x+1)^5.$ Since 
$\Delta\left(2,3,8\right)$ is maximal 
the group of automorphisms $G'$ has
order 96 (see {\bf Example \ref{group96}}).  This leads to case 
$B.3$ in {\bf Theorem \ref{maintheorem}}.

{\underline{Subcase 4.3}}: We now consider the index 6 extension
$$
\begin{array}{cccccc}
\Delta (n,4n,4n)&\lrar&\bbz_{4n}\\
 \downarrow&&\downarrow\\
\Delta (2,3,4n)&\lrar&G'
\end{array}
$$
According to (\cite{BC}, T8) we must have $z_1\neq z_2$ in order
to be able to extend further.  Moreover $n=2,3$ or $6$ and 
$\{z_1,z_2,z_3\}= \{T,T^4,T^{-5}\}$. Up to a permutation we can take
$a=1$, $[4n,b]=1, [4n,c]=4$. In fact $c=4$
since the only element of order $n$ in $\{T,T^4,
T^{-5}\}$ is $T^4.$  
\\
$\bullet$ If $n=2$ then
$b=3$ and the surface is $C_1: y^8 = x(x-1)^3(x+1)^4$ with automorphism 
group $G'$ of order $48$.  This leads to case $E.1$ in
{\bf Theorem \ref{maintheorem}}.
\\
$\bullet$ If $n=3$ then $b=7$ and the surface is
$C_2: y^{12} = x(x-1)^4(x+1)^7$ with automorphism group $G'$ of order $72$.
This leads to case $E.2$ in {\bf Theorem \ref{maintheorem}}.
\\
$\bullet$ If $n=6$ then $b=19$ and the  surface is 
$C_3: y^{24}=x(x-1)^4(x+1)^{19}$, with automorphism group $G'$ of order $144$.
This leads to case $E.3$ in {\bf Theorem \ref{maintheorem}}.
The automorphism group of $ y^{24}= x(x-1)^4(x+1)^{19}$ is given as
$central(\bbz_6):S_4$ (also described in \cite{BC} with more on this in \cite{KSS}).  On
the other hand, it is easy to see that $C_2=C_3/\bbz_2$ and that
$C_1=C_3/\bbz_3$ (as quotient surfaces). Since the cyclic group actions
are central, the group structures of Aut$(C_1)$ and Aut$(C_2)$ follow
directly.

{\underline{Subcase 4.4}}: Finally we look at the extension 
described in row 12 of the GS table (index 4)
$$
\begin{array}{cccccc}
\Delta (n,2n,2n)&\lrar&\bbz_{2n}\\
 \downarrow&&\downarrow\\
\Delta (2,4,2n)&\lrar&G'
\end{array}
$$
As usual, $a=1$, $[2n,b]=1$ and $[2n,c]=2$. $\Delta\left(2,4,2n\right)$
is maximal  except when $n=4$, but this was treated in subcase 4.2.  Thus
assume $n\neq 4.$

Two possibilities arise: either $z_1=z_2$ (and an
extension always exists), or $z_1\neq z_2$ and $4|n$ (see \cite{BC},
T9).  \\
$\bullet$ Suppose $z_1=z_2$ and hence $a=b=1$, $c=2n-2$.  Notice that
in this case the surface can be brought to the form
$y^{2n}=x^{2n-2}(x-1)(x+1)$, and by moving the branch point at $x=0$
to $\infty$,
it has the affine equation $y^{2n}= x^2-1$ (``Accola-Maclachlan'' type).
This leads to case $A.2$ in  {\bf Theorem \ref{maintheorem}}.\\
$\bullet$ Suppose $z_1\neq z_2$ and hence $b\neq 1$. Then there
is $k$ such that $k^2\equiv 1\ (mod\ 2n)$ and 
$k\not\equiv 1\ (mod\ 2n)$.  We can take $k=b$ and thus $c=2n-1-k.$   
From the fact that $[2n,c]=2$ it follows that $[2n,k+1]=2.$  But then
the congruences above imply that $k=n+1.$
The surface in this case is then $y^{2n} =
x(x-1)^{n+1}(x+1)^{n-2}$ (known as the Kulkarni surface).  
More on this surface in {\bf Example
\ref{kulkarni}}.  This leads to case $B.2$ of
{\bf Theorem \ref{maintheorem}}.


\vskip 10pt
{\bf\Large\S5 Full Automorphism Groups of Lefschetz
Surfaces}

Based on the calculations in the previous section, we now
determine the automorphism groups of the generalized Lefschetz surfaces 
$$L: y^n = x^a(x-1)^b(x+1)^c,\ 1\le a,b,c\le n-1,\ a+b+c\equiv 0\ (mod\ n)$$ 

First of all, we need to address the question of maximality of the
extensions worked out in \S4. In that section, we analyzed when a skep
$\theta:\Delta\lrar\bbz_n$ extends to $\Delta_1\lrar G_1$ where the
extension $\Delta\subset\Delta_1$ is taken from row 1 of the GS table
(case 3), row 13 (case 5), row 3 (case 4.1), row 7 (case 4.2), row 11 (case
4.3) and row 12 (case 4.4). We now must determine if $\Delta_1\lrar
G_1$ further extends to $\Delta_2\lrar G_2$. Of course this can
happen only if $\Delta_1$ is not maximal in the GS table.

For example, the $\bbz_7$ action on the surface $C:
y^7=x(x-1)^2(x+1)^4$ is uniformized by $\theta: \Delta
(7,7,7)\lrar\bbz_7$.  According to $\underline{\hbox{Case}\ 3}$ (\S4),
since $1+2+2^2\equiv 0\ (mod\ 7)$, there is an index 3 extension (coming
from row 1 of the GS table), $\theta_1:\Delta (3,3,7)\lrar G_1$, where
$G_1$ is a group of order $21=3\times 7$ acting on the surface.  It
turns out that $C$ (of genus $3$) is the Klein surface with 168
automorphisms (see {\bf Lemma \ref{klein2}}). The extra index 8 extension 
arises from row 6 in the GS table. What happens in this case
is that we have consecutive extensions
$$\Delta (7,7,7)\fract{r_1}{\lrar} (3,3,7)\fract{r_6}{\lrar}
\Delta (2,3,7)$$ 
where $r_i$ refers to row $i$ of the GS table. The resulting extension
corresponds to $r_4$ of the GS table. The following is
well-known.

\ble \label{klein2}
The only cyclic covering of $\bbp^1$ of genus 3 with an automorphism
group of order $168$ is the curve $y^7 = x(x-1)^2(x+1)^4$ . It is
isomorphic to Klein's surface $x^3y + y^3 + x= 0$, with automorphism
group $PSL(2,7).$  
\ele

{\sc Proof}: An extension $\Delta (2,3,7)\lrar G$ of $\Delta
(7,7,7)\lrar\bbz_7$ occurs if and only if $\{z_1, z_2,z_3\}=\{T,
T^2,T^4\}$ (or an equivalent triple) for some generator $T$ of
$\bbz_7$ (\cite{BC}, case T1).  The group $G$ has order $168$ in this
case and is isomorphic to PSL$(2,7)$. In fact, this shows  that
if $y^n = x^a(x-1)^b(x+1)^c$ has genus $3$ and an automorphism group
of order $168$, then necessarily $n=7$, $a=1, b=2, c=4$ (up to
permutation).  Finally, by a change of variables (see \cite{K}) the surface 
$K: x^3y + y^3 + x = 0$ is equivalent to
to $C': y^7 = x(x-1)^2$, which in turn is
equivalent to $C: y^7 = x(x-1)^2(x+1)^4$.\hfill\za

The complete list of consecutive extensions from the GS table is:

\begin{enumerate}
\item 
$\Delta (n,n,n)\fract{r_1}{\lrar}\Delta (3,3,n)\fract{r_3}{\lrar}
\Delta (2,3,2n)$, $n\ge 4$, is equivalent to $r_2.$
\vspace{0.05in}
\item 
$\Delta (7,7,7)\fract{r_1}{\lrar}\Delta
(3,3,7) \fract{r_6}{\lrar}\Delta (2,3,7)$
is equivalent to $r_4.$
\vspace{0.05in}
\item 
$\Delta (9,9,9)\fract{r_{1}}{\lrar}\Delta (3,3,9)
\fract{r_{13}}{\lrar}\Delta (2,3,9)$ is equivalent to
to $r_{9}$.
\vspace{0.05in}
\item 
$\Delta (n,2n,2n)\fract{r_{3}}{\lrar}\Delta (2,2n,2n)
\fract{r_{3}}{\lrar}\Delta (2,4,2n)$, $n\ge 3$, is equivalent to
to $r_{12}$.
\vspace{0.05in}
\item 
$\Delta (4,8,8)\fract{r_{3}}{\lrar}\Delta (2,8,8)
\fract{r_{11}}{\lrar}\Delta (2,3,8)$ is equivalent to
to $r_{7}$.
\vspace{0.05in}
\item 
$\Delta (n,n,n)\fract{r_{3}}{\lrar}\Delta (2,n,2n)
\fract{r_{14}}{\lrar}\Delta (2,3,2n)$, $n\ge 4$, is equivalent to
to $r_{2}$.
\vspace{0.05in}
\item $\Delta(n,4n,4n)\fract{r_{3}}{\lrar}\Delta(2,2n,4n)
\fract{r_{14}}{\lrar}\Delta(2,3,4n)$, $n\ge 2$, is equivalent to $r_{11}.$
\vspace{0.05in}
\item 
$\Delta (4,8,8)\fract{r_{12}}{\lrar}\Delta (2,4,8)
\fract{r_{14}}{\lrar}\Delta (2,3,8)$ is equivalent to
to $r_{7}$.
\end{enumerate}

A skep $\Delta (n,n,n)\lrar\bbz_n$ does not extend to $\Delta (2,3,2n
)\lrar G$, according to \cite{BC}, N7, and so cases (1) and (6) can be
dismissed.  The second case above is covered by {\bf Lemma \ref{klein2}}.
Cases (5) and (8) are equivalent and are covered by subcase (4.2).
Calculations in \cite{BC} imply that case (3) not arise.
That leaves cases (4) and (7), which are covered by subcases
(4.4) and (4.3) resp. 

The above discussion together with the calculations in the previous
section completely determines the automorphism groups of the cyclic
covers $y^n=x^a(x-1)^b(x+1)^c$ and is summarized in
{\bf Theorem \ref{maintheorem}}. It remains to discuss the structure of the
groups in question and this we do in the following list of examples.

\bex\label{periodthree} In this example we discuss the action of 
$\bbz_n\ltimes\bbz_3$ on $C:
y^{n} = x(x-1)^k(x+1)^{k^2}$, where $1+k+k^2\equiv 0\ (mod\ n)$ (and hence
$k^3\equiv 1\ (mod\ n )$). The following argument is adapted from (\cite{L},
p. 177). First we consider the {\sl equivalent} curve
$$y^{n} = (x-1)(x-j)^{k}(x-j^2)^{k^2}$$
where $j=e^\frac{2\pi i}{ 3}$.  Set $1+k+k^2=\alpha n$ and
$k^3 = \beta n + 1$. The $\bbz_3$ action is given by
$$S:(x,y)\mapsto 
\left(jx,j^{\alpha}y^{k}(x-j^2)^{-\beta }\right).
$$
This is well defined (i.e. it does act on $C$ and has period $3$).
Let $T$ be the order
$n$-transformation. Then $ST = T^kS$.  The group generated by $S$ and
$T$ is a non-abelian semi-direct product $\bbz_n\ltimes\bbz_3$.
This is case $C_1$ of 
{\bf Theorem \ref{maintheorem}}.
\eex

\bex \label{accola-maclachlan} 
In $\underline{\hbox{case 4.4}}$ of \S4, the 
surface $C: y^{2n}=x^{2n-2}(x-1)(x+1)$ 
was found to have an additional $\bbz_2$
action (besides the involution $(x,y)\mapsto (-x,y)$).
We describe an element of order 4 acting on this surface. First, $C$ is 
isomorphic to the surface $C':y^{2n}=x^2-1$
Consider the automorphism
$$\displaystyle
u: (x,y)\mapsto \left(\frac{x}{y^n}, \frac{\zeta}{ y}\right),\ \ \zeta^n=-1$$
It is easy to see that $u^2(x,y)=(-x,y)$ and hence $u^4=1$. Let
$v$ be the $\bbz_{2n}$ cyclic generator $(x,y)\lrar (x,\zeta y)$.
The full automorphism of this curve is now given by
$$<u,v\ |\ u^4=v^{2n}=(uv)^2=[u^2,v]=1 >. $$
The subgroup generated by $u^2$ is central and moding out by it we get
the dihedral group $D_{4n}=<x, v | x^2=v^{2n}= (xv)^2=1>$.  This is
case $A_2$ of {\bf Theorem \ref{maintheorem}}.
See also {\bf Example \ref{special7}}.
\eex

\bex\label{twistedz2} In $\underline{\hbox{Case 4.1}}$, \S4, we
established the existence of a $\bbz_2$ action on curves of the form
$C: y^n = x(x-1)^b(x+1)^{n-b-1}$ with $b^2\equiv 1\ (mod\ n)$. We now
explicitly describe such an action. First of all, $C$ is isomorphic to
$C': y^n = (x+1)^b(x-1)$. Consider the
automorphism
\begin{eqnarray*}\displaystyle
&&u: (x,y)\mapsto (-x, (-1)^ly^b(x+1)^{-\beta}),\ \mbox{where}\\
&& b^2 = \beta n + 1,
nl\equiv b+1\ (mod\ 2)\ \mbox{and}\ l+bl\equiv \beta\ (mod\ 2)
\end{eqnarray*}
It can be checked that $u$ is a well defined involution, provided
$n\not\equiv 0\ (mod\ 8).$  One can also check that $uvu=v^b,$
 where $v$ is the usual cyclic action $(x,y)\mapsto (x,\zeta y)$.
This calculation in fact describes the action of the semi-direct product 
 $\bbz_n\ltimes\bbz_2=<u,v\ | u^2=v^n=1, uvu=v^b>$  on the surface.
This is case $B_1$ of {\bf Theorem  \ref{maintheorem}}.
\eex

\bex \label{kulkarni}
The Kulkarni surface $y^{2n}=(x+1)^{n+1}(x-1),$ discussed
in $\underline{\hbox{case 4.4}}$ (\S4), has an automorphism group of
order $8n = 8g+8$ presented by  (see \cite{BC}, T9)
$$\left<u,v\ |\ u^4=v^{2n}=(uv)^2=u^2vu^2v^{n-1}=1\right>\approx
\left(\mathbb Z_{2n}\ltimes\mathbb Z_2\right):Z_2$$
The subgroup generated by $\{u^2, v\}$ is the semi-direct product
$\mathbb Z_{2n}\ltimes\mathbb Z_2$. 
\eex

\bex\label{group96} The genus $3$ surface
$C:y^8=x(x-1)^2(x+1)^5$ discussed in $\underline{\hbox{case
4.2}}$ (case B.3) has a group of automorphisms which we claim is a split
extension of $\bbz_4\oplus\bbz_4$ by $S_3$. This group is a
permutation group on 12 letters, generated by the cycles
\begin{eqnarray*}\displaystyle & &(1,4)(2,7)(3,10)(5,8)(6,11)(9,12)\\ 
& &(1,10,9,5)(2,4,11,3,7,12,6,8)\\
& &(1,2,3)(4,5,6)(7,8,9) (10,11,12)
\end{eqnarray*} 
See \cite{BC}, T4. Running {\sc magma} on
this group shows that there is one single normal subgroup $H$ of order
$16$ with quotient a non-abelian group of order $6$ (necessarily
$S_3$). $H$ is abelian with $7$ subgroups of order $4$. It is then
necessarily $\bbz_4\oplus\bbz_4$ as claimed
.\eex

{\bf Classical Lefschetz surfaces}: Now let $C$ be given as a prime
Galois cover of the sphere. Such curves are necessarily given by the
equation $y^p = f(x)$ (see \cite{K}) and if $C$ is belyi it is isomorphic
to the curve
$$y^p = x^{a}(x-1)^{b}(x+1)^{c},~~~a+b+c\equiv 0\ (mod\ p), ~~
1\leq a,b,c \le p-1$$
Note that $p$ being prime is crucial in deriving this form of the
equation.  For example $\bbz_{12}$ acts on the curve $x^3 + y^4=1$
with quotient $\bbp^1$ and its equation is not amenable (via
birational transformations) to the form above. In fact by moving one
of the branch points to $\infty\in\bbp^1$ we can get the following
biholomorphic model

\ble \label{model2}
A Lefschetz surface is birationally equivalent to one with
equation
$$y^p = x^a(x+1),~1\leq a < \frac{p-1}{2}$$  
\ele
{\sc Proof}: By {\bf Lemma \ref{model}}, $C:y^p = x^a(x-1)^b(x+1)^c$ is
isomorphic to $y^p =x^{a'}(x-1)(x+1)^{c'}$ (to see this raise to a power
$l$ such that $lb\equiv 1\ (mod\ p)$ and then reduce $a,c$  modulo $p$).  
Note that $a'+ c' + 1=p$,
and so by the change of variables $X=x^{-1},~Y=\eta yx^{-1}$,
$\eta^p=-1$, we obtain the isomorphic curve $Y^p = (X-1)(X+1)^{c'}$.
This is equivalent to the curve $y^p=x^{c'}(x+1)$. 

So without loss of generality, a Lefschetz curve has the form $C_a: y^p
= x^a(x+1)$ for some $a$, $1\leq a\le p-1$.  Applying the change of
variables $Y=yx^{-1}, X=x^{-1}$, we find that $C_a$ is isomorphic to
$C_{p-a-1},$ and hence one can choose $1\leq a\leq \frac{p-1}{2}$.
It remains to see that $y^p=x(x+1)$ is isomorphic to $y^p=x^\frac{p-1}{
2}(x+1)$.  According to {\bf Lemma \ref{model}}, $y^p=x^{p-2}(x+1)(x-1)$ is
isomorphic to $y^p=x(x-1)^\frac{p-1}{ 2}(x+1)^\frac{p-1}{ 2}$ (since
$\left(\frac{p-1}{ 2}\right)(p-2)\equiv 1\ (mod\ p))$, and hence to
$y^p=x^\frac{p-1}{2}(x-1)^\frac{p-1}{ 2}(x+1)$. Another change of variables
leads to the equation $y=(x-1)^{\frac{p-1}{2}}(x+1)$, which is
what we wanted to prove.  \hfill\za

\noindent{\bf Remark 3}:
\label{unique}
Notice that according to the above representation there is a {\sl
unique} Lefschetz surface when $p=5$ (namely $y^5=x(x+1)$). In fact it
is the unique surface of genus $\frac{p-1}{ 2}=2$ which admits a
$\bbz_5$ action! See \cite{RR}.  Indeed if $C$ has genus 2 and $\bbz_5$
acts on it, then $C/\bbz_5\cong\bbp^1$ and so $C:y^5=f(x)$. It has
exactly 3 branch points and so must be Lefschetz, and hence of
the form above. Its group of automorphism is in fact $\bbz_{10}$ (see
below).  

\bth
Let $C$ be a Lefschetz surface isomorphic to $y^p = x^a(x+1),$
$1\leq a< \frac{p-1}{2}$. Let $G=$Aut$(C)$. Then\\
(1) If $a=1$, $G$ is cyclic of order $2p$.\\
(2) If $p\equiv 1\ (mod\ 3)$, $p>7$ and $1+a+a^2\equiv 0\ (mod\ p)$, 
then $G$ is the unique non-abelian group of order $3p$.\\
(3) If $p=7$ and $a=2$, then $G=PSL(2,7)$ is the simple group of order
168.\\
(4)  For all other cases, $\aut (C)=\bbz_p$.
\end{theorem}

{\sc Proof}: We simply read off the possibilities from the
classification table in {\bf Theorem \ref{maintheorem}} when $n=p$ is a prime
(only cases $A.1, B.1, C.1, C.2$ apply). In that table,
representative curves have equations of the form $C: y^p =
x(x-1)^b(x+1)^c$, or $y^p=x^a(x+1)$, where 
$ac\equiv 1\ (mod\ p).$  Without loss of
generality we can assume $1\le a<(p-1)/2.$
 We can rule out case $B_1$ right away as it would
imply $b=p-1,c=0.$   \\ 
$A.1$ Here $a=1$ and $Aut\ C\approx \mathbb Z_{2p}.$ Indeed 
$C: y^p=x(x+1)$ has the obvious
$\bbz_2$ action, $(x,y)\mapsto (-x-1,y),$ commuting
with the $\mathbb Z_p$ action.\\ 
C.1 Here $p\equiv 1\ (mod\ 3)$, $p>7$ and $1+a+a^2\equiv 0\ (mod\ p)$.  
For all such curves $\bbz_3$ acts in a
twisted fashion (see {\bf Example \ref{periodthree}}).\\ 
C.2 Here $p=7$, $a=2$ and
the surface is Klein's curve $y^7=x^2(x+1)$ ({\bf Example \ref{klein}}).  
\hfill\za

\bex \label{klein} The surface $y^7=x^2(x+1)$ is isomorphic to Klein's
curve $x^3y + y^3z + z^3x = 0$, written in projective coordinates (see 
{\bf Lemma
\ref{klein2}}).  The action of $\bbz_3$ is given by permuting $x,y$ and
$z$ in a $3$-cycle. An automorphism of order 7 is given by
$x\mapsto\zeta x, y\mapsto\zeta^4y$ and $z\mapsto\zeta^2z$, where
$\zeta$ is a primitive $7^{th}$ root of unity. The action of $\bbz_2$ is a
bit more involved and is described for instance in \cite{Ba}.  \eex

\noindent{\bf Remark 4}:
The representation of a Lefschetz curve as $y^p=x^a(x+1)$ with
$1\leq a<\frac{p-1}{2}$ is generally not unique. When $p=7$ it is true
however that there only two isomorphism classes of Lefschetz curves:
the Klein surface $y^7=x^2(x+1)$ and the surface $y^7=x(x+1)$.  

\bpr \label{conditions}
Two Lefschetz surfaces $y^p=x^a(x+1)$ and $y^p=x^b(x+1)$, 
$1\le a,b <\frac{p-1}{ 2},$ are isomorphic if and only if one of the following
is true: $a=b$, $ab+b+1\equiv 0\ (mod\ p)$, $ab+a+1\equiv 0\ (mod\ p)$,
$a+b+ab\equiv 0\ (mod\ p)$ or 
$ab\equiv 1\ (mod\ p)$.
\epr

The proof is an easy check using {\bf Lemma \ref{model}}, whose converse is
true in the prime case.  Note that the
count of isomorphism classes of Lefschetz surfaces is given in \cite{RR}
with an extension for more general $p$-elliptic curves in \cite{KSS}.


\vskip 10pt
{\bf\Large\S6. Curves of Fermat Type}

These are curves of the form
$${ F}: y^n + x^d = 1,\ \ 1< d\leq n$$
The affine curves are smooth.
When $d=n$ the projective curve is smooth as well and has genus
$\frac{(n-1)(n-2)}{ 2}$. However, when $d<n$ there is branching over
$\infty$, and a little calculation shows that in this case 
$g = \frac{1}{2}(2-d-[d,n] + (d-1)n)$.  
Notice that there is an obvious action of
$\bbz_d\oplus\bbz_n$ on $ F$, with each
cyclic action having quotient the Riemann sphere.

Fermat curves are a special class of cyclic covers of
the Riemann sphere, but 
our prior techniques for the study of Aut$({ F})$ don't apply directly
when $d\geq 4$ (because the map 
$F\to F/\mathbb Z_n$ has branching over more than 
$3$ points).  One way to get around this problem is to consider
instead the quotient map
$$\beta: { F}\fract{\bbz_d\oplus\bbz_n}{\ra 4}\bbp^1$$
This turns out to be a belyi map. To see this
consider the $\bbz_n$ quotient given by the map $(x,y)\mapsto x$.
This is a regular $n$-fold covering branched over  
$B =\{x\ |\
x^d=1\}$ and also over $x=\infty$ when $d<n$.  The map
$f:\bbp^1\lrar\bbp^1, x\mapsto x^d$, is branched over $0, \infty$ and
maps all of $B$ to 1. Therefore the composite $\beta$
$$\label{bcc}
\begin{array}{cccccc}
\beta: &F & \stackrel{\psi}{\lrar} &\bbp^1 &\stackrel{f}{\lrar}
&\bbp^1\\
 &(x,y)& \lrar&x & \lrar &x^d
\end{array}
$$
is branched only over $0,1$ and $\infty$ and hence is belyi. As can easily
be seen this composite is exactly the quotient map ${ F}\lrar
{ F}/\bbz_d\oplus\bbz_n$. Observe that
$|\beta^{-1}(0)|=n,|\beta^{-1}(1)|=d$. The branching
over $\infty$ is different.  By the change of variable
$z=1/x$, we can bring this equation (locally at $x=\infty,\ z=0$ ) to 
the form $v^n = z^{n-d}(z^d-1)=u^{m}$, where $v=zy$ and 
$m=n-d$ (\cite{Mi}, p. 73). 
As explained in \S3, one needs
$[m,n]=[d,n]$ ramification points over $\infty$ for the $\bbz_n$ quotient
$\psi: { F}\lrar\bbp$ and just one ramification point for $f$.
This means that over $\infty$ there are $[d,n]$ ramification points
for $\beta$. Therefore, we can
uniformize the action of $G=\bbz_d\oplus\bbz_n$ on $F$ by the 
triangle group 
$\displaystyle
\Delta\left (\frac{nd}{n},\frac{nd}{ d},\frac{nd}{ [d,n]}\right)=
\Delta (d,n,(n,d)).
$

In summary, the action of $\bbz_d\oplus\bbz_n$ on { F} is uniformized
by an exact sequence
\begin{eqnarray}\label{ses2}
1\to\Pi\to\Delta (d,n, (n,d))\stackrel{\theta}
{\rightarrow}\bbz_d\oplus\bbz_n\to 1,
\end{eqnarray}
where $\Pi\cong \pi_1({ F})$.
To determine $\aut ({ F})$ (for all $d$ and $n$), we need to study the
existence and extendability of skeps $\Delta (d,
n, (n,d))\fract{\theta}{\lrar}\bbz_d\oplus\bbz_n$.  This in part has
been carried out in a recent preprint \cite{BCC}.

In {\bf Examples \ref{special7}} and {\bf\ref{special8}} we consider 
the cases $d=2$ and $3$, but for now we assume $d\ge 4$. 
If \mbox{$d\not{|}\ n$} then $1<d<n<(d,n)$ and an examination of the GS table 
reveals that the skep 
$\displaystyle
\Delta(d,n,(d,n))\stackrel{\theta}{\to}\mathbb Z_d\oplus\mathbb Z_n$
is not extendable.  The only relevant cases are 13 and 14, but these 
correspond to $d=3$ and $d=2$ resp.  Therefore, if $d\!\not{|}\ n$ then
$\aut(F)=\mathbb Z_d\oplus\mathbb Z_n.$
  
Thus we are seeking extensions of the skep 
$\displaystyle
1\to\Pi\to\Delta(d,n,n)\to \mathbb Z_d\oplus\mathbb Z_n\to 1
$
for those cases where $d|n, d\ge 4.$
This case-by-case study has been carried out in \cite{BCC}.  It turns out
that  extensions 
occur for cases (1), (2) and (3) of the GS table, and no others. 
We tacitly assume $4\le d\le n$ and $d|n$ in what 
follows. 

\ble (Extension E.1) \label{e1}
Any skep $\Delta (d, n,n)\fract{\theta}{\lrar}\bbz_d\oplus\bbz_n$ admits
an index 2 extension  $\Delta (2,n,2d)\rightarrow G^{\prime}$. If we set 
$\theta(x_1)=S,\theta(x_2)=T$ then $G^{\prime}$ has the presentation
$$<S,T,U\ |\ S^d=T^n=1, ST=TS,U^2=1,USU=S,UTU=S^{-1}T^{-1}>\approx
\left(\mathbb Z_d\oplus\mathbb Z_n\right)\ltimes\mathbb Z_2$$ 
  \ele

\noindent{\sc Proof}:
The abelianization of $\Delta (d,n,n)$ is $\bbz_d\oplus\bbz_n$ and $\theta$
is the abelianization homomorphism up to an automorphism of 
$\bbz_d\oplus\bbz_n.$ 
The extension we are seeking has the form
$$\begin{array}{ccccccccc}
1&\lrar&\Pi&\lrar &\Delta(d,n,n)&\fract{\theta}{\lrar}&\mathbb Z_d\oplus\mathbb Z_n &\lrar& 1\\
&&\decdnar{=}&&\decdnar{}&&\decdnar{}\\
1&\lrar&\Pi&\lrar &\Delta(2,n,2d)&\fract{\theta^{\prime}}{\lrar}&G^{\prime}&\lrar& 1
\end{array}
$$
We can choose elliptic generators $y_1,y_2,y_3$  of $\Delta (2,n,2d)$ 
 so that
$$x_1=y_1^2,x_2=y_2,x_3=y_3y_2y_3, y_3^2=1,y_1y_2y_3=1$$ 
Set $U=\theta^{\prime}(y_3).$ If the above diagram exists then
compatibility requires $USU=S,UTU=S^{-1}T^{-1}$, so $G^{\prime}$ has
the stated presentation.  Conversely, if $G^{\prime}$ has the given
presentation then the extension exists.

\hfill\za

\ble (Extension E.2, see \cite{BCC})\label{bujacirreconder}
Any skep $\Delta (n,n,n)\lrar G=\bbz_n\oplus\bbz_n$ extends to $\Delta
(3,3,n)$.
\ele

Making use of these calculations we can prove.

\bth\label{fermat2}
Let $F$ be the surface given by $y^n+x^d=1$, where $4\leq d \leq n$,
and let Aut$(F)$ be its group of automorphisms.\\
(1) If $d$ does not divide $n$, then Aut$(F)=\bbz_d\oplus\bbz_n$. \\
(2) If $d=n$, Aut$(F)$ is the semi-direct product of 
$\bbz_n\oplus\bbz_n$ with the symmetric group $S_3$.  \\ 
(3) If $d|n$, $d<n$, Aut$(F)$ is the semi-direct product of
$\bbz_d\oplus\bbz_n$ by $\bbz_2$. A presentation is
$$ \hbox{Aut} (F) = 
\big < S,T,U\ \bigm |\ S^d=T^n=1,ST=TS,U^2=1, USU=S, UTU=S^{-1}T^{-1}\big >$$
\end{theorem}

{\sc Proof}:

If $d=n$ then the conditions for E.1 and E.2 are satisfied and an extension 
of index 6 always occurs (see {\bf Remark 5}). Indeed $S_3$ acts on 
$x^n+y^n=1$ (by permuting $x,y$ and $z$ in the equivalent projective equation
$x^n+y^n+z^n=0$) and so it follows that $(\bbz_n\oplus\bbz_n)\ltimes S_3$
acts on $F$ as the full automorphism group.
The only other case to consider is the one where $d|n,1<d<n,$ but this is
covered by {\bf Lemma \ref{e1}}.
 \hfill\za

\noindent{\bf Remark 5}:
\label{degree6}
The index 6 extension of $\Delta (n,n,n)\to\bbz_n\oplus\bbz_n$ to
$\Delta (2,3,2n)\to (\bbz_n\oplus\bbz_n)\ltimes S_3 $ is a combination 
of E.1 and E.2.

\bex\label{special7}
In this example we consider the curve $F:y^n+x^2=1.$ By counting 
ramification orders we see that the genus is given by
$$\displaystyle
g=\left\{\begin{array}{ll}
			(n-2)/2 &\mbox{if $n$ is even}\\
			(n-1)/2 &\mbox{if $n$ is odd}
	\end{array}
  \right.
$$

The action of $\mathbb Z_2\oplus\mathbb Z_n$ is uniformized by
$\displaystyle 
\Delta(2,n,2n)\to\mathbb Z_2\oplus\mathbb Z_n\cong\mathbb Z_{2n}$
if $n$ is odd and by
 $\displaystyle 
\Delta(2,n,n)\to\mathbb Z_2\oplus\mathbb Z_n$
if $n$ is even.  In the odd case the only possible extension comes from case 14 in the GS table, but according to \cite{BC} no such extension is possible.

Thus assume $n$ is even.  Then, there are 2 cases to consider in the GS table, namely cases 3 and 11 (for $ n=8$).  The possible extensions of the skep
$\displaystyle 
\Delta(2,n,n)\to\mathbb Z_2\oplus\mathbb Z_n
$
are to $\Delta(2,4,n)\to G^{\prime}$ and  $\Delta(2,3,8)\to G^{\prime}$
respectively.  The last case is impossible (see \cite{BCC}) and the first case
always exists ({\bf Lemma \ref{e1}}).  In this case 
$G^{\prime}\cong (\mathbb Z_2\oplus\mathbb Z_n)\ltimes\mathbb Z_2$ has the presentation
$$G^{\prime}=
\left<a,b,u\ |\ a^n=1,b^n=1,(ab)^2=1,ab=ba,u^2=1,uau=b,ubu=a\right>$$
\eex
This agrees with  part 3 of {\bf Theorem \ref{fermat2} }. 
\bex\label{special8}
Now consider the curve $F:y^n+x^3=1$ of genus
$$g=\left\{\begin{array}{ll}
			n-2 &\mbox{if $n\equiv 0\ (mod\ 3)$ }\\
			n-1 &\mbox{otherwise}
	\end{array}
  \right.
$$
The $\mathbb Z_3\oplus\mathbb Z_n$ action on $F$ is uniformized by
$\displaystyle
\Delta(3,n,n)\to\mathbb Z_3\oplus\mathbb Z_n 
$ if $n\equiv 0\ (mod\ 3)$ and by $
\Delta(3,n,3n)\to\mathbb Z_3\oplus\mathbb Z_n \cong \mathbb Z_{3n}
$ if $n\not{\!\equiv} 0\ (mod\ 3).$

If $n\equiv 0\ (mod\ 3)$ the possible extensions arise from cases 3, 11 and 12 
of the GS table.  Cases 11 and 12 do not extend (see \cite{BCC}).  Case 3 
concerns the possibility of extending 
$\displaystyle
\Delta(3,n,n)\to\mathbb Z_3\oplus\mathbb Z_n$ to 
$\Delta(2,6,n)\to G^{\prime}.$  This is possible, and in fact 
$G^{\prime}\cong(\mathbb Z_3\oplus\mathbb Z_n)\ltimes\mathbb Z_2$ has
the presentation
$$\displaystyle \left<a,b,u\ |\ a^n=b^n=(ab)^3=1, ab=ba,u^2=1,
uau=b,ubu=a\right>$$
This follows from {\bf Lemma \ref{e1}}.

Now assume  $n\not{\!\!\equiv} 0\ (mod\ 3).$  The only possible
extension of 
$\Delta(3,n,3n)\to\mathbb Z_3\oplus\mathbb Z_n \cong \mathbb Z_{3n}$
comes from case 13 of the GS table: namely an extension to
$\Delta(2,3,3n)\to G^{\prime}.$ This extension exists if, and only if,
$n=4.$  In this case $G^{\prime}$ is a central extension of $\mathbb Z_4$ 
by the
alternating group $A_4$.  A presentation for $G^{\prime}$ is
$$\displaystyle
G^{\prime}\cong\left<u_1,u_2\ |\ u_1^2=u_2^3=[u_1,(u_1u_2)^3]=1 \right>$$
The central subgroup $\mathbb Z_4$ is generated by $(u_1u_2)^3.$  See
\cite{BC} for the details.
\eex

\noindent{\sc Acknowledgements}: The first author would like to thank
the University of British Columbia and the PIms Institute for their
hospitality while this work was being conducted.

\small
\addcontentsline{toc}{section}{Bibliography}
\bibliography{biblio}
\bibliographystyle{plain[8pt]}

\vskip 30pt
\flushleft{Sadok Kallel \\
Universit\'e des Sciences et Technologies de Lille\\
U.F.R de Math\'ematiques\\
59655 Villeneuve d'Ascq, France\\
{\sc Email:}  sadok.kallel@agat.univ-lille1.fr}

\vskip 10pt
\flushleft{
Denis Sjerve\\
Department of Mathematics\\
University of British Columbia\\
Vancouver V6T 1Z2\\
{\sc Email:}  sjer@math.ubc.ca}

\end{document}